\documentclass[12pt]{amsart}
\usepackage[mathscr]{eucal}
\usepackage{amscd}
\usepackage{amsfonts}
\usepackage{amsmath, amsthm, amssymb}
\usepackage{latexsym}
\usepackage[dvips]{graphics}
\usepackage{enumerate}
\usepackage{cite}
\begin{document}
\title[Relation-theoretic contraction principle]
{Relation-theoretic contraction principle in metric-like as well as partial metric spaces}
\author[Ahmadullah, Khan and Imdad ]{ Md Ahmadullah$^{1}$, Abdur Rauf Khan$^{2}$ and Mohammad Imdad$^{3}$}
%\thanks{$^{\ast}$Corresponding author}
\maketitle
\begin{center}
{\footnotesize $^{1}$Department of Mathematics, Aligarh Muslim University, Aligarh-202002, U.P., India.\\
E-mail addresses: {\tt ahmadullah2201@gmail.com; mahmadullah.rs@amu.ac.in}\\
$^2$Department of Mathematics, College of Science, Jazan University, Jazan, Kingdom of Saudi Arabia and
   University of Polytechnics, Applied Mathematics Section, Aligarh Muslim University, Aligarh-202002, U.P., India.
   E-mail addresses: {\tt abdurrauf@rediffmail.com}\\
$^{3}$Department of Mathematics, Aligarh Muslim University, Aligarh-202002, U.P., India.\\
E-mail addresses: {\tt mhimdad@yahoo.co.in; mhimdad@amu.ac.in}}
\end{center}

\begin{abstract} % abstract
In this paper, we extend the Banach contraction principle to metric-like as well as partial metric spaces (not essentially complete)
equipped with an arbitrary binary relation. Thereafter, we derive some fixed point results which are sharper versions of
the corresponding known results of the existing literature. Finally, we use an example to demonstrate our main result.
\end{abstract}

\vspace{.3cm}\noindent

\vspace{.3cm}\noindent

\section{Introduction}
Fixed point theory is a relatively old but still a young research area amongst all non-linear mathematical sciences.
Banach contraction principle remains a celebrated result in metric fixed point theory which was established by
Banach \cite{bnch1922} in 1922. In recent years, many researchers studied fixed point results in ordered metric spaces
(e.g., \cite{AhmadJI,Alamimdad,Alamimdad2,NietoL2005,RanR2004,rus2008,Turinicid1986} and references cited therein).
The much discussed idea of metric space has been generalized and improved by introducing several variants such as:
metric-like space, partial metric space, symmetric space, pseudo metric space, b-metric space, 2-metric space, G-metric space and several others.

\vspace{0.3cm} In 1994, Matthews \cite{matt1994} introduced the concept of partial metric space and
also extended Banach contraction principle in such spaces. In recent years, a multitude of metrical
fixed point theorems were extended to the partial metric (e.g., \cite{altun2011,Ilic2011,kara2011,matt1994,oltra2004,rus2008}
and references cited therein) and such research activity is still on.

\vspace{0.3cm} Hitzler \cite{hitz2001}, proved an interesting extension of the Banach contraction principle by introducing
dislocated metric spaces. Here it can be pointed out that dislocated metric spaces are also referred as metric-like spaces
(e.g., Amini-Harandi \cite{amini2012}). For further details on metric-like spaces
one can consult \cite{amini2012,aydi2015,chen2015,hitz2001,hitz2000,huseyin2013,kara2013} and references cited therein.

\vspace{0.3cm}
The aim of this paper is to extend the Banach contraction principle to metric-like spaces (not essentially complete)
equipped with an arbitrary binary relation. Thereafter, we derive some fixed point results which are sharper versions of
the corresponding known results of the existing literature. Finally, we furnish an example to demonstrate our main result.

\vspace{0.3cm} Throughout this paper, $\mathbb{R}^+$, $\mathbb{N}$ and $\mathbb{N}_{0}$ respectively, stand for the set of non-negative real numbers,
the set of natural numbers and the set of whole numbers.

\section{Preliminaries}
We begin with definitions of partial metric and metric-like spaces well followed by some of their relevant properties.

\vspace{.3cm} \noindent
{\bf Definition 1.} \cite{matt1994} Let $X$ be a non-empty set. Then a mapping $p : X \times X \to \mathbb{R}^{+}$ is said to be
	partial metric on $X$ if for all $x, y, z \in X$,
	\begin{enumerate}
		\item [{${(p_1)}$}] $x = y \iff p(x, x) = p(x, y) = p(y, y),$
		\item [{${(p_2)}$}]  $p(x, x) \leq p(x, y),$
		\item [{${(p_3)}$}] $p(x, y) = p(y, x),$
		\item [{${(p_4)}$}] $p(x, y) \leq p(x, z) + p(z, y) - p(z, z).$
	\end{enumerate}

The pair $(X,p)$ is called a partial metric space.

\vspace{.3cm} \noindent
{\bf Definition 2.} \cite{hitz2001} Let $X$ be a non-empty set. Then a mapping $\sigma : X \times X \to \mathbb{R}^{+}$ is said to be
	metric-like (or dislocated) on $X$ if for all $x, y, z \in X$
	\begin{enumerate}
		\item [{${(\sigma_1)}$}] $\sigma(x, y)=0 \Rightarrow x = y,$
		\item [{${(\sigma_2)}$}] $\sigma(x, y) = \sigma(y, x),$
		\item [{${(\sigma_3)}$}] $\sigma(x, y) \leq \sigma(x, z) + \sigma(z, y).$
	\end{enumerate}

The pair $(X,\sigma)$ is called a metric-like (or dislocated) space. Here it can be pointed out that all the requirements of a metric are met out except $\sigma(x, x)$ may be positive for $x\in X$.

\vspace{.3cm} \noindent
{\bf Remark 1.}\label{rmk2.1}
	Obviously every metric is partial metric and every partial metric is metric-like but converse implication is not true in general.

\vspace{.3cm} \noindent
{\bf Example 1.}\label{exp2.1}
	Let $X=\{a,b,c\}$ and $\sigma, p : X\times X \to \mathbb{R}^+$ given by
	$$\begin{cases}
	\sigma(a,a)=\sigma(b,b)=0,\\
	\sigma(c,c)=\sigma(a,b)=\sigma(b,a)=2,\\
	\sigma(a,c)=\sigma(c,a)=\sigma(b,c)=\sigma(c,b)=1;
	\end{cases}$$
and $$p(x,y)=\begin{cases}
	0, & \hbox{{$x=y=a$};}\\
	1, & \hbox{otherwise.}
	\end{cases}$$
	
	\vspace{0.3cm} Then $(X,\sigma)$ is a metric-like space but not a partial metric space due to the fact that $\sigma(c,c)=2\nleq 1=\sigma(c,a)$
	while $(X,p)$ a partial metric space but not a metric space as $p(b,b)\ne 0.$

\vspace{.3cm}
We utilize the following terminologies which are essentially available in Amini-Harandi \cite{amini2012}.

\vspace{.3cm} \noindent
{\bf Definition 3.} \cite{amini2012}
	Let $(X,\sigma)$ be a metric-like space and $\{x_n\}$ a sequence in $X$. Then we say that
	\begin{itemize}
		\item the sequence $\{x_n\}$ converges to a point $x$ in $X$ if and only if $\displaystyle\lim_{n\to \infty} \sigma(x_n, x)=\sigma(x, x)$,
		\item the sequence $\{x_n\}$ is Cauchy in $X$ if and only if $\displaystyle\lim_{n,m\to \infty} \sigma(x_n, x_m)$ is exist and finite,
		\item the metric-like space $(X,\sigma)$ is complete if every Cauchy sequence $\{x_n\}$ in $X$ converges to a point $x$ in $X$ with respect to $\tau_\sigma$ (topology generated by $\sigma$) such that
		$$\lim_{n,m\to \infty} \sigma(x_n, x_m)=\sigma(x, x)=\lim_{n\to \infty} \sigma(x_n, x).$$
	\end{itemize}

\vspace{0.3cm}
Next, we present some relevant relation-theoretic notions:

\vspace{.3cm} \noindent
{\bf Definition 4.} \cite{Lips1964} Let $X$ be a non-empty set $X$. Then a subset of $X\times X$ is called a binary relation on $X$, which will be denoted by $\mathcal{R}$. We say that ``$x$ is related to $y$ under $\mathcal{R}$" if and only if $(x,y)\in \mathcal{R}$.

\vspace{.3cm}
In what follows, $\mathcal{R}$ stands for a non-empty binary relation.

\vspace{.3cm} \noindent
{\bf Definition 5.}
	\cite{Maddux2006} A binary relation $\mathcal{R}$ on $X$ is called complete if for all $x,y$ in $X$, either $(x,y)\in
	\mathcal{R}$ or $(y,x)\in \mathcal{R}$ which is denoted as $[x,y]\in \mathcal{R}$.

\vspace{.3cm} \noindent
{\bf Definition 6.}
	\cite{Alamimdad} Let $f$ be a self-mapping defined on a non-empty
	set $X$. Then a binary relation $\mathcal{R}$ on $X$ is called $f$-closed if
	$(fx,fy)\in \mathcal{R}~\text{whenever}~(x,y)\in \mathcal{R},~ {\rm for~ all}~ x,y\in X.$

\vspace{.3cm} \noindent
{\bf Definition 7.}
	\cite{Alamimdad} Let $\mathcal{R}$ be a binary relation defined on a non-empty set $X$. Then a sequence $\{x_n\}$ in $X$ is called $\mathcal{R}$-preserving if
	$(x_n,x_{n+1})\in\mathcal{R},\;\;{\rm for~all}~n\in \mathbb{N}.$

\vspace{.3cm}
Motivated by Alam and Imdad \cite{Alamimdad2}, we introduce relation-theoretic variants of completeness and continuity in metric-like spaces.

\vspace{.3cm} \noindent{\bf Definition 8.}
	Let $(X,\sigma)$ be a metric-like space and $\mathcal{R}$ a binary relation on $X$.
	We say that $(X,\sigma)$ is $\mathcal{R}$-complete if every $\mathcal{R}$-preserving Cauchy sequence $\{x_n\}$ in $X$,
	there is some $x \in X$ such that
	$$\lim_{n,m\to \infty} \sigma(x_n, x_m)=\sigma(x, x)=\lim_{n\to \infty} \sigma(x_n, x).$$

Recall that the limit of a convergent sequence in metric-like spaces is not necessarily unique.

\vspace{.3cm} \noindent
{\bf Remark 2}\label{rmk2.2}. Every complete metric-like space is
	$\mathcal{R}$-complete but not conversely.
	The notion of $\mathcal{R}$-completeness coincides with completeness if the relation $\mathcal{R}$ is universal.

\vspace{.3cm} \noindent
{\bf Definition 9.}
	Let $(X,\sigma)$ be a metric-like space. Then a mapping
	$f:X\rightarrow X$ is said to be continuous-like at $x$ if $f(x_n)\stackrel{\tau_\sigma}{\longrightarrow} f(x)$ for
	any sequence $\{x_n\}$ with $x_n\stackrel{\tau_\sigma}{\longrightarrow} x$.
	As usual, $f$ is said to be continuous-like if it is continuous-like on the whole $X$.

\vspace{.3cm} \noindent
{\bf Definition 10.}
	Let $(X,\sigma)$ be a metric-like space and $\mathcal{R}$ a binary relation on $X$. Then a mapping
	$f:X\rightarrow X$ is said to be $\mathcal{R}$-continuous-like at $x$ if $f(x_n)\stackrel{\tau_\sigma}{\longrightarrow} f(x)$ for
	any $\mathcal{R}$-preserving sequence $\{x_n\}$ with
	$x_n\stackrel{\tau_\sigma}{\longrightarrow} x$.
	As usual, $f$ is said to be $\mathcal{R}$-continuous-like if it is $\mathcal{R}$-continuous-like on the whole $X$.

\vspace{.3cm} \noindent
{\bf Remark 3.}\label{rmk2.3} On metric-like spaces, every continuous mapping is continuous-like and every continuous-like mapping is
	$\mathcal{R}$-continuous-like but not conversely. The notion of
	$\mathcal{R}$-continuity-like coincides with continuity-like if the relation $\mathcal{R}$ is universal.

\vspace{.3cm} \noindent
{\bf Definition 11.} \cite{Alamimdad}
	Let $(X,\sigma)$ be a metric-like space and $\mathcal{R}$ a binary relation on $X$. Then $\mathcal{R}$ is said to be
	$\sigma$-self-closed if for any $\mathcal{R}$-preserving sequence
	$\{x_n\}$ with $x_n\stackrel{\tau_\sigma}{\longrightarrow} x$, there
	is a subsequence $\{x_{n_k}\}$ of $\{x_n\}$ such that $[x_{n_k},x]\in\mathcal{R},~{\text for ~all}~k\in \mathbb{N}.$

\vspace{.3cm} \noindent
{\bf Definition 12.}
	\cite{SametT2012} Let $(X,\sigma)$ be a metric-like space and $\mathcal{R}$ a binary relation
	on $X$. Then a subset $D$ of $X$ is said to be $\mathcal{R}$-directed if for every pair of points $x,y \in D$, there
	is $z$ in $X$ such that $(x,z)\in\mathcal{R}$ and $(y,z)\in\mathcal{R}$.

\vspace{.3cm} \noindent
{\bf Definition 13.}
	\cite{KBR2000} Let $(X,\sigma)$ be a metric-like space, $\mathcal{R}$ a binary relation on $X$ and $x,y$ a pair of points in $X$. Then a finite sequence $\{z_0,z_1,z_2,...,z_{l}\}$ in $X$ is said to be a path
	of length $l$ (where $l\in \mathbb{N}$) joining
	$x$ to $y$ in $\mathcal{R}$ if $z_0=x, z_l=y$ and $(z_i,z_{i+1})\in\mathcal{R}$ for each $i\in \{1,2,3,\cdots ,l-1\}.$

\indent Here it can be pointed out that a path of length $l$ involves $(l+1)$ elements of $X$
that need not be distinct in general.

\vspace{.3cm}
\indent In a metric-like space $(X,\sigma)$, a self-mapping $f$ on $X$ and a binary relation $\mathcal{R}$ on $X$, we employ the following notations:

\begin{itemize}
	\item $F(f)$: the set of all fixed points of $f$;
	\item $\Upsilon(x,y,\mathcal{R})$: the family of all paths joining $x$ to $y$ in $\mathcal{R}$.
\end{itemize}

\section{Main results}
\vspace{.3cm} \noindent
 {\bf Theorem 1.}\label{th2.1} {\it
		Let $(X,\sigma)$ be a metric-like space equipped with a binary relation $\mathcal{R}$ defined on $X$ and $f$ a self-mapping on
		$X$. Suppose that the following conditions are satisfied:
		\begin{enumerate}
			\item [$(a)$] there exists a subset $Y\subseteq X$ with $fX\subseteq Y$ such that $(Y,\sigma)$ is $\mathcal{R}$-complete,
			\item [$(b)$] there exists $x_0$ such that $(x_0, fx_0)\in \mathcal{R}$,
			\item [$(c)$] $\mathcal{R}$ is $f$-closed,
			\item [$(d)$] either $f$ is $\mathcal{R}$-continuous-like or $\mathcal{R}|_Y$ is $\sigma$-self-closed,
			\item [$(e)$] there exists a constant $k\in [0,1)$ such that $(\text{for all} ~x,y\in X\;\textrm{with}\; (x,y)\in \mathcal{R})$
			$$\sigma(fx,fy)\leq k \sigma(x,y).$$
		\end{enumerate}
		Then $f$ has a fixed point. Moreover, if
		\begin{enumerate}
			\item [${(f)}$] $\Upsilon(fx,fy,\mathcal{R}^s)$~ is non-empty, for each $x,y\in X$.
		\end{enumerate}
		Then $f$ has a unique fixed point. }
\begin{proof} Construct Picard iterate $\{x_n\}$ corresponding to $x_0$, $i.e, x_n=f^nx_0$\; for all $n\in \mathbb{N}_0$. Since $(x_0,fx_0)\in
	\mathcal{R}$ and $\mathcal{R}$ is $f$-closed, we find that
	$$(fx_0,f^2x_0),(f^2x_0,f^3x_0),\cdots ,(f^nx_0,f^{n+1}x_0),\cdots\in
	\mathcal{R},$$ so that
	\begin{equation}\label{eq2.1}
	(x_n,x_{n+1})\in \mathcal{R}\;\;{\rm for~all}~n\in \mathbb{N}_0.
	\end{equation}
	Hence $\{x_n\}$ is an
	$\mathcal{R}$-preserving sequence. Our assert that $\{x_n\}$ is Cauchy sequence. To establish this, using the condition $(e)$, we have (for all $n\in \mathbb{N}_0$)
	$$\sigma(x_{n+1},x_{n+2})=\sigma(fx_{n},fx_{n+1})\leq k\sigma(x_{n},x_{n+1}),$$
	which yields by induction that
	\begin{equation}\label{eq2.2}
	\sigma(x_{n+1},x_{n+2})\leq k^{n+1} \sigma(x_0,fx_0)\;\forall~n\in \mathbb{N}_0.
	\end{equation}
	By triangular inequality, (\ref{eq2.2}) and for all $n,m\in \mathbb{N}_0$ with $m>n$, we have
	\begin{eqnarray*}
		\nonumber \sigma(x_{n},x_{m})&\leq& \sigma(x_{n},x_{n+1})+\sigma(x_{n+1},x_{n+2})+\cdots+\sigma(x_{m-1},x_{m})\\
		&\leq&(k^{n}+k^{n+1}+\cdots+k^{m-1})\sigma(x_0,fx_0)\\
		&=& k^{n}\sigma(x_0,fx_0)\sum\limits_{j=0}^{m-n-1} k^{j}\\
		&\leq& \frac{k^{n}}{1-k}\sigma(x_0,fx_0)\\
		&\rightarrow& 0\;{\rm as}\; n\rightarrow \infty,
	\end{eqnarray*}
	which shows that the sequence $\{x_n\}$ is Cauchy in $Y$. Hence, the sequence $\{x_n\}$ is an $\mathcal{R}$-preserving Cauchy in $Y$.
	By $\mathcal{R}$-completeness of $(Y,\sigma)$, there exists $y\in Y$ such that the sequence $\{x_n\}$ converges
	to $y$ with respect to $\tau_\sigma$ (topology generated by $\sigma$) $i.e.,$
	\begin{equation}\label{eq2.3}
	\displaystyle\lim_{n \to \infty}\sigma(x_{n},y)=\sigma(y,y)=\displaystyle\lim_{n,m \to \infty}\sigma(x_{n},x_{m})=0.
	\end{equation}
	Firstly, assume that $f$ is $\mathcal{R}$-continuous-like. Then ~~$x_{n+1}=fx_n\stackrel{\tau_\sigma}{\longrightarrow} fy$, so that
	\begin{eqnarray}\label{eq2.4}
	\lim_{n\to \infty}\sigma(x_{n+1},fy)=\lim_{n\to \infty}\sigma(fx_{n},fy)=\sigma(fy,fy)=\lim_{n,m \to \infty}\sigma(x_{n},x_{m})=0.
	\end{eqnarray}
	On using triangular inequality, (\ref{eq2.3}) and (\ref{eq2.4}), we have
	$\sigma(y, fy)=0,$
	so that $y$ is a fixed point of $f$.

	\noindent Alternately, if $\mathcal{R}|_Y$ is $\sigma$-self-closed, then due to the fact that
	$\{x_n\}$ is an $\mathcal{R}$-preserving sequence in $Y$ and
	$x_n\stackrel{\tau_{\sigma}}{\longrightarrow} y, \exists$ a subsequence
	$\{x_{n_k}\}{\rm \;of\;} \{x_n\} \;{\rm with}\;\;[x_{n_k},y]\in\mathcal{R},~~~\text{for~all} ~k\in \mathbb{N} _0.$
	%Notice that, $[x_{n_k},y]\in\mathcal{R}\;~~\forall~k\in \mathbb{N} _0$ implies that either
	%$(x_{n_k},y)\in\mathcal{R}\;\;~\forall~k\in \mathbb{N} _0$
	%or, $(y, x_{n_k}) \in\mathcal{R}\;\;~~~\forall~k\in \mathbb{N} _0$.
	In view of the condition $(e)$ and the symmetry of the metric-like $\sigma$, we have
	$$\sigma(x_{n_{k+1}}, fy)=\sigma(fx_{n_{k}}, fy)\leq k\sigma(x_{n_{k}}, y).$$
	Taking the limit as $k\to \infty$ and using (\ref{eq2.3}), we have
	\begin{equation}\label{eq2.5}
	\displaystyle\lim_{k\to \infty}\sigma(x_{n_{k+1}}, fy)=0.
	\end{equation}
	On using triangular inequality, (\ref{eq2.3}) and (\ref{eq2.5}), we have
	$$\sigma(y, fy)=0,$$
	so that, $fy=y, ~i.e.,~ y$ is a fixed point of $f$.
	
	\vspace{0.3cm} Next, if $F(f)$ is singleton then result follows.
	Otherwise, take $p,q$ be two arbitrary elements of $F(f)~  i.e.,$
	$$fp=p\;{\rm and}\;fq=q.$$
	Owing to the condition $(f)$, there exists a path
	(say $\{x_0,x_1,x_2, \cdots, x_l\}$) of some finite length $l$ in
	$\mathcal{R}^s$ from $p$ to $q$ such that
	$$x_0=p,\;x_{l}=q\;{\rm and}\;[x_i,x_{i+1}]\in\mathcal{R}\;{\rm for~each}\;i\in\{1,2,\cdots ,l-1\}.$$
	As $\mathcal{R}$ is $f$-closed, we have
	$$[f^nx_i,f^nx_{i+1}]\in\mathcal{R}\;{\rm for~each}\;i\in\{1,2,\cdots ,l-1\}\;{\rm and}\;{\rm for~each}\;n\in
	\mathbb{N}_0.$$
	Now, on using triangular inequality and hypothesis $(e)$, we obtain
	\begin{eqnarray*}
		\nonumber \sigma(p,q)&=& \sigma(f^nx_{0},f^nx_{l})\\
		&\leq&\sum\limits_{i=0}^{l-1} \sigma(f^nx_i,f^nx_{i+1})\\
		&\leq& k\sum\limits_{i=0}^{l-1} \sigma(f^{n-1}x_i,f^{n-1}x_{i+1})\\
		&\leq&k^2\sum\limits_{i=0}^{l-1}\sigma(f^{n-2}x_i,f^{n-2}x_{i+1})\\
		&&\vdots \\
		&\leq& k^n\sum\limits_{i=0}^{l-1}\sigma(x_i,x_{i+1})\\
		&\rightarrow& 0\;{\rm as}\; n\rightarrow \infty,
	\end{eqnarray*}
	yielding thereby $p=q$. Hence $f$ has a unique fixed point.
\end{proof}

Particularly, by setting $Y=X$ in Theorem 1, we deduce the following:

\vspace{.3cm}\noindent
{\bf Corollary 1.}\label{cor2.1} {\it
		Let $(X,\sigma)$ be a metric-like space equipped with a binary relation $\mathcal{R}$ on $X$ and $f$ a self-mapping on
		$X$. Suppose that the following conditions are satisfied:
		\begin{enumerate}
			\item [$(a)$] $(X,\sigma)$ is $\mathcal{R}$-complete,
			\item [$(b)$] there exists $x_0$ such that $(x_0, fx_0)\in \mathcal{R}$,
			\item [$(c)$] $\mathcal{R}$ is $f$-closed,
			\item [$(d)$] either $f$ is $\mathcal{R}$-continuous-like or $\mathcal{R}$ is $\sigma$-self-closed,
			\item [$(e)$] there exists a constant $k\in [0,1)$ such that $(\text{for~ all}~ x,y\in X\;\textrm{with}\; (x,y)\in \mathcal{R})$
			$$\sigma(fx,fy)\leq k \sigma(x,y).$$
		\end{enumerate}
		Then $f$ has a fixed point. Moreover, if
		\begin{enumerate}
			\item [${(f)}$] $\Upsilon(fx,fy,\mathcal{R}^s)$~ is non-empty, for each $x,y\in X$.
		\end{enumerate}
		Then $f$ has a unique fixed point. }

\vspace{0.3cm} In view of Remarks 2 and 3, we deduce the following natural result:

\vspace{.3cm} \noindent
{\bf Corollary 2.}\label{cor2} {\it
	Let $(X,\sigma)$ be a metric-like space equipped with a binary relation $\mathcal{R}$ defined on $X$ and $f$ a self-mapping on
	$X$. Suppose that the following conditions are satisfied:
	\begin{enumerate}
		\item [$(a)$] there exists a subset $Y\subseteq X$ with $fX\subseteq Y$ such that $(Y,\sigma)$ is complete,
		\item [$(b)$] there exists $x_0$ such that $(x_0, fx_0)\in \mathcal{R}$,
		\item [$(c)$] $\mathcal{R}$ is $f$-closed,
		\item [$(d)$] either $f$ is continuous or $\mathcal{R}|_Y$ is $\sigma$-self-closed,
		\item [$(e)$] there exists a constant $k\in [0,1)$ such that $(\text{for all} ~x,y\in X\;\textrm{with}\; (x,y)\in \mathcal{R})$
		$$\sigma(fx,fy)\leq k \sigma(x,y).$$
	\end{enumerate}
	Then $f$ has a fixed point. Moreover, if
	\begin{enumerate}
		\item [${(f)}$] $\Upsilon(fx,fy,\mathcal{R}^s)$~ is non-empty, for each $x,y\in X$.
	\end{enumerate}
	Then $f$ has a unique fixed point. }
	
\vspace{0.3cm} In respect of ``$\mathcal{R}^s$-directedness of $fX$'' and ``completeness of $\mathcal{R}|_{fX}$'', we mention the following corollary.

\vspace{0.3cm} \noindent
{\bf Corollary 3.}\label{cor2.2} {\it
		Theorem 1 remains true if we replace the condition $(f)$ by one of the following
		conditions besides retaining the rest of the hypotheses:
		\begin{enumerate}
			\item [${(f^\prime)}$] $fX$ is $\mathcal{R}^s$-directed,
			\item [${(f^{\prime\prime})}$] $\mathcal{R}|_{fX}$ is complete.
		\end{enumerate}}
	\begin{proof}
		Assume that ${(f^\prime)}$ holds. Then for each pair of points $a, b$ in $fX$, $\exists x \in X$ such that $[a,x]\in\mathcal{R}$ and
		$[b,x]\in\mathcal{R}$ so that the finite sequence $\{a,x,c\}$ is a path of length 2 joining $a$ to $b$ in $\mathcal{R}^s$.
		Thus, for each $a,b \in fX$, $\Upsilon(a,b,\mathcal{R}^s)$ is non-empty and hence the conclusion follows from Theorem 1.
		
		Now, let us consider condition $(f^{\prime\prime})$ holds. Then for each pair of points $a, b \in fX$,
		$[a,b]\in\mathcal{R}$, which implies that $\{a,b\}$ is a path
		of length 1 from $a$ to $b$ in $\mathcal{R}^s$,
		so that for each $a,b \in fX$, $\Upsilon(a,b,\mathcal{R}^s)$ is non-empty.
		Finally, on the proceeding lines of Theorem 1, the conclusion is immediate.
	\end{proof}
	
	Under the universal relation $\mathcal{R}$, Corollary 2 gives rise an improved version of Banach contraction principle on the metric-like space, as follows:
	
	\vspace{0.3cm} \noindent
	{\bf Corollary 4.}\label{cor2.3} {\it
			Let $(X,\sigma)$ be a metric-like space and $f$ a self-mapping on
			$X$. Suppose that there exists a subset $Y$ of $X$ with $fX\subseteq Y\subseteq X$ such that $(Y,\sigma)$ is complete. If
			there is $k\in [0,1)$ such that
			$$\sigma(fx,fy)\leq k \sigma(x,y),~{\text for ~all}~ x,y\in X,$$
			then $f$ has a unique fixed point.}
	
\vspace{0.3cm}\noindent	{\bf Remark 4.}
		Corollary 4 is an improved version of  Corollary 3.6 due to Aydi and Karapinar \cite{aydi2015}.

\vspace{0.3cm}\noindent	{\bf Remark 5.}
		In view of Remark \ref{rmk2.1}, the class of metric-like spaces remains relatively larger than classes of partial metric spaces and metric spaces.
		Consequently, on can easily deduce the analogues results corresponding to Theorems 1 and
		Corollaries  1-4 in the settings of partial metric spaces and metric spaces.
	
	\vspace{0.3cm} Let $\Omega$ be the set of all mappings $\rho : \mathbb{R}^{+} \to \mathbb{R}^{+}$ such that
	\begin{enumerate}
		\item [${(i)}$] $\rho$ is a Lebesgue-integrable on each compact subset of $\mathbb{R}^{+}$, and
		\item [${(ii)}$] $\int_{0}^{\varepsilon}\rho(t)>0, ~\text{for~all}~\varepsilon>0.$
	\end{enumerate}
	
	\vspace{0.3cm} Now, we predict the following generalized form of Theorem 1 employing the integral type contractive condition, which runs as follows.
	
	\vspace{0.3cm}\noindent {\bf Theorem 2.}\label{thm2.3} {\it
			Let $(X,\sigma)$ be a metric-like space equipped with a binary relation $\mathcal{R}$ defined on $X$ and $f$ a self-mapping on
			$X$. Suppose that:
			\begin{enumerate}
				\item [$(a)$] there exists a subset $Y\subseteq X$ with $fX\subseteq Y$ such that $(Y,\sigma)$ is $\mathcal{R}$-complete,
				\item [$(b)$] there exists $x_0$ such that $(x_0, fx_0)\in \mathcal{R}$,
				\item [$(c)$] $\mathcal{R}$ is $f$-closed,
				\item [$(d)$] either $f$ is $\mathcal{R}$-continuous-like or $\mathcal{R}|_Y$ is $\sigma$-self-closed,
				\item [$(e)$] there exists $\rho \in \Omega$ such that $(\text{for all}~ x,y\in X\;\textrm{with}\; (x,y)\in \mathcal{R})$
				$$\int_{0}^{\sigma(fx,fy)}\rho(t)dt\leq k\int_{0}^{\sigma(x,y)}\rho(t)dt.$$
				
			\end{enumerate}
			Then $f$ has a fixed point. Moreover, if
			\begin{enumerate}
				\item [${(f)}$] $\Upsilon(fx,fy,\mathcal{R}^s)$~ is non-empty, for each $x,y\in X$.
			\end{enumerate}
			Then $f$ has a unique fixed point. }

	\begin{proof}
		Proceeding on the lines of Theorem 1 and Theorem 2.1 of Branciari \cite{branc2002}, on can complete the proof. Hence the full proof is not included.
	\end{proof}
	
Under universal relation, Theorem 2 reduces to an improved version of Theorem 2.1 due to Branciari \cite{branc2002}, which runs as follows
	
\vspace{0.3cm}\noindent	{\bf Corollary 5.}\label{cor2.4} {\it
			Let $(X,\sigma)$ be a metric-like space and $f$ a self-mapping on
			$X$. Suppose that there exists a subset $Y$ of $X$ with $fX\subseteq Y\subseteq X$ such that $(Y,\sigma)$ is complete. If
			there exists $k\in [0,1)$ and  $\rho \in \Omega$ such that (for all $x,y \in X$)
			$$\int_{0}^{\sigma(fx,fy)}\rho(t)dt\leq k\int_{0}^{\sigma(x,y)}\rho(t)dt,$$
			then $f$ has a unique fixed point.}

\section{Illustrative Example}
	
\vspace{0.3cm} Finally, we furnish an example demonstrating the usability of Theorem 1.
	
	\vspace{0.3cm}\noindent	{\bf Example 2.}
		Let $X=\{a, b, c\}$ and $\sigma : X\times X \to \mathbb{R}^+$ defined by
		$$\begin{cases}
		\sigma(a,a)=\sigma(b,b)=0, ~\sigma(c,c)=3,\\
		\sigma(a,b)=\sigma(b,a)=1,\\
		\sigma(a,c)=\sigma(c,a)=\sigma(b,c)=\sigma(c,b)=2.
		\end{cases}$$
		Then $(X,\sigma)$ is a metric-like space which is neither a metric nor a partial metric space.
		Now, we define a mapping $f: X\to X$ by
		$$fa=b, ~~fb=b~~ \text{and}~fc=a,$$ a binary relation $\mathcal{R}=\{(a,a),(b,b),(a,b)\}$ on $X$ and $Y=\{a,b\}$.
		Then $\mathcal{R}$ is $f$-closed and $Y$ is $\mathcal{R}$-complete. Take any $\mathcal{R}$-preserving sequence
		$\{x_n\}$ with $$x_n\stackrel{\tau_\sigma}{\longrightarrow} x ~{\rm and}~
		(x_n,x_{n+1})\in\mathcal{R} ~{\rm for~ all}~ n\in \mathbb{N} _0.$$
		Notice that if $(x_n,x_{n+1})\in\mathcal{R}|_Y$ for all $n\in\mathbb{N} _0,$ then there exists an integer $N \in\mathbb{N} _0$
		such that $x_n=p \in \{a,b\} ~\forall n\geq N$. So, we can take a subsequence
		$\{x_{n_k}\}$ of the sequence $\{x_n\}$ with $x_{n_k}=p ~\forall k\in \mathbb{N} _0$, which
		amounts to saying that $[x_{n_k},p]\in \mathcal{R}|_Y~ \forall k\in \mathbb{N} _0$. Therefore,
		$\mathcal{R}|_Y$ is $\sigma$-self-closed.
		
		\vspace{0.3cm} In order to check the condition $(e)$ of Theorem 1, it is sufficient to show that the condition $(e)$ holds for $x\in\{a,b\}$ and $y=c$ (or, $x=c$ and $y\in\{a,b\}$), as in rest of the cases $\sigma(fx,fy)=0$.
		If $x\in\{a,b\}$ and $y=c$, then $\sigma(fx,fy)=1\leq k2= k\sigma(x,y)$ holds for all $k\in[\frac{1}{2},1)$.
		As $\mathcal{R}|_{fX}$ is complete, the condition $(f)$ holds. Thus all the requirements of Theorem 1 are fulfilled.
		Notice that $f$ has a unique fixed point (namely, `$b$\text{'}).

	\vspace{0.3cm}\noindent
	{\bf Competing interests}\\
	The authors declare that they have no competing interests.
		
	\vspace{0.3cm}\noindent {\bf Authors' contributions}\\
	The authors contributed equally in this paper.

\end{document}